\documentclass[11pt]{article}
\usepackage{amsfonts}
\usepackage{latexsym}
\usepackage{amsmath}
\usepackage{amssymb}
\usepackage{amsthm}
\newtheorem{theorem}{Theorem}[section]
\newtheorem{lemma}[theorem]{Lemma}

\newtheorem{proposition}[theorem]{Proposition}

\theoremstyle{definition}

\newtheorem{remark}{Remark}
\theoremstyle{remark}

\newtheorem*{note*}{Note}

\makeatletter
\newcommand{\ls}{\leqslant}
\newcommand{\gr}{\geqslant}

\makeatother

\begin{document}

\small

\title{\bf On the approximation of a polytope by its dual $L_{p}$-centroid bodies
}

\author{Grigoris\ Paouris\thanks{Partially supported by an NSF grant (DMS-0906150) } \   and Elisabeth M. Werner 
\thanks{Partially supported by an NSF grant, a FRG-NSF grant and  a BSF grant}}

\date{}

\maketitle

\begin{abstract}
\footnotesize We show that the rate of convergence on the approximation of volumes of a convex symmetric polytope $P\in \mathbb R^n$ by its dual $L_{p}$-centroid bodies is independent of the geometry of $P$. In particular we show that if $P$ has volume $1$,
$$ \lim_{p\rightarrow \infty} \frac{p}{\log{p}} \left( \frac{|Z_{p}^{\circ}(P)|}{|P^{\circ}|} -1 \right) = n^{2} . $$
We provide an application to the approximation of polytopes by uniformly convex sets.
\end{abstract}

\bigskip

\section{Introduction}

\noindent 
Let $K$ be a convex body in $\mathbb R^n$ of volume $1$ and,  for $\delta\in (0,1)$,  let  $K_{\delta}$ be  the convex floating body of $K$ \cite{SW2}. It is  the intersection of all halfspaces $H^{+}$ whose defining hyperplanes $H$ cut off a set of volume $\delta$ from $K$. Note that $K_{\delta}$ converges to $ K$ in the Hausdorff metric as $\delta\rightarrow 0$. C. Sch\"utt and the second name author showed an exact formula for the convergence of   volumes \cite{SW2},
$$ \lim_{\delta\rightarrow 0} \frac{ |K|- |K_{\delta}|}{ \delta^{\frac{2}{n+1}}} = {\rm as}_{1}(K), $$
which involves  the affine surface area of $K$, ${\rm as}_{1}(K)$. The same phenomenon (and similar formulas) has been observed for other types of approximation using instead of floating bodies, convolution bodies \cite{Schm}, illumination bodies \cite{W1} or Santal\'o bodies \cite{MW}. We refer to e.g. \cite{GaZ}, \cite{GrZh}-\cite{Klain2}, \cite{Lud}-\cite{LZ}, \cite{SW1}-\cite{SA2}, \cite{WY}-\cite{Z3} for  further details, extensions and applications.  Another family of bodies that  approximate a given convex body $K$ are the $L_{p}$-centroid bodies of $K$ introduced by Lutwak and Zhang \cite{LZ}. For a symmetric convex body $K$ of volume $1$ in $\mathbb R^n$ and $1\ls p \ls n$, the $L_{p}$-centroid body $Z_{p}(K)$ is the convex body that has support function
$$ h_{Z_{p}(K)}(\theta) = \left(\int_{K} |\langle x, \theta\rangle|^{p} d x \right)^{\frac{1}{p}} , \ \theta \in S^{n-1}. $$
Note that $Z_{p}(K)$ converges  to $K$ in the Hausdorff metric as $p\rightarrow \infty$. It has been shown in \cite{PW} that the family of $L_{p}$-centroid bodies is isomorphic to the family of the floating bodies:  $K_{\delta}$ is isomorphic to $Z_{\log{\frac{1}{\delta}}}(K)$. However, it was proved in \cite{PW} that in the case of $C^{2}_{+}$ bodies, the convergence of volume of the $L_{p}$-centroid bodies is independent of the ``geometry" of $K$: For any symmetric convex body in $\mathbb R^n$ of volume $1$ that is $C_{+}^{2}$ (i.e. $K$ has $C^{2}$ boundary with everywhere strictly positive Gaussian curvature),
$$ \lim_{p\rightarrow\infty} \frac{p}{\log{p}} \left(  |Z_{p}^{\circ}(K)| - |K^{\circ}|\right) = \frac{n(n+1)}{2} |K^{\circ}| . $$
In this work we show that the same phenomenon occurs also in the case of polytopes. We show the following
\begin{theorem} \label{theorem1}
\noindent Let $K$ be a symmetric polytope of volume $1$ in $\mathbb R^n$. Then 
$$  \lim_{p\rightarrow\infty} \frac{p}{\log{p}} \left(  |Z_{p}^{\circ}(K)| - |K^{\circ}|\right) = n^{2} |K^{\circ}| . $$
\end{theorem}
\noindent As an application of this result we get bounds for the approximation of a polytope by a uniformly convex body with respect to the symmetric difference metric:
\begin{theorem} \label{theorem2}
\noindent  Let $ P$ be a symmetric  polytope in $\mathbb R^n$. Then there exists $p_{0}=p_{0}( P) $ such that for every $p\gr p_{0}$, there exists a $p$-uniformly convex body $K_{p}$ such that 
$$  d_{s}(P, K_{p}) \ls 2n^{2} |P| \frac{\log{p}}{p},$$
where $d_{s}$ is the symmetric difference metric. 
\end{theorem}

\smallskip
\noindent
The statements and proofs are for symmetric convex bodies only. 
If $K$ is not symmetric, then $Z_{p}(K)$ does not converge to $K$ since the $Z_{p}(K)$ are centrally symmetric by definition.
However, all results  can be extended to the non-symmetric case with minor modifications of the proofs 
by using the non-symmetric  version of the $L_p$-centroid bodies from \cite{Lud} (see also \cite{HabSch}).
\vskip 3mm
\noindent 
The paper is organized as follows. In section 2 we give some bounds for the approximation of volume in the case of a general convex body. In section 3 we consider the case of polytopes and we give the proof of Theorem \ref{theorem1}. Finally,  in  section 4, we discuss  approximation of a polytope by $p$-uniformly convex bodies  (see \cite{LindTz}) and we give the proof of Theorem \ref{theorem2}. 

\vskip 4mm

\noindent 
{\bf Notation.}
\vskip 2mm
\noindent
We work in ${\mathbb R}^n$, which is
equipped with a Euclidean structure $\langle\cdot ,\cdot\rangle $.
We denote by $\|\cdot \|_2$ the corresponding Euclidean norm, and
write $B_2^n$ for the Euclidean unit ball and $S^{n-1}$ for  the
unit sphere. Volume is denoted by $|\cdot |$. 
We write $\sigma $ for the rotationally
invariant surface measure on $S^{n-1}$. 
\par
\noindent 
A convex body is a compact convex subset $C$ of ${\mathbb R}^n$ with
non-empty interior. We say that $C$ is symmetric, if $x\in
C$ implies that $ -x\in C$. We say that $C$ has center of mass at the
origin if $\int_C\langle x,\theta\rangle dx=0$ for every $\theta\in
S^{n-1}$. The support function $h_C:{\mathbb R}^n\rightarrow
{\mathbb R}$ of $C$ is defined by $h_C(x )=\max\{\langle x,y\rangle
:y\in C\}$. 
$C^{\circ}  = \{ y\in {\mathbb R}^n: \langle x,y\rangle \ls 1\;\hbox{for all}\; x\in C\} 
$ is the polar body of $C$.
\par
\noindent 
We refer to \cite{Gar} and \cite{Sch} for basic facts from the
Brunn-Minkowski theory.

\vskip 3mm

\noindent \textbf{Acknowledgments.} The authors would like to thank the American Institute of Mathematics. Part of this work has been carried out during a stay at AIM.

 \vskip 3mm

\section{General Bounds}

\noindent Let $K$ be  a symmetric convex body in $\mathbb R^n$ of volume $1$. Let $\theta \in S^{n-1}$. We define the parallel section function $f_{K,\theta} :[-h_{k}(\theta), h_{k}(\theta)] \rightarrow \mathbb R_{+}$ by
$$ f_{K,\theta} (t) := | K \cap ( \theta^{\perp} + t \theta) | . $$

\noindent 
By Brunn's principle,  $ f_{K,\theta} ^{\frac{1}{n-1}}$ is concave and attains its maximum at $0$. So we have that
\begin{equation}\label{Brunn}
 \left(1- \frac{t}{h_{K}(\theta)} \right)^{n-1}  f_{K,\theta}(0) \ls f_{K,\theta}(t) \ls f_{K,\theta}(0) .  
 \end{equation}

\noindent 
The right-hand side inequality is sharp if and only if $K$ is a cylinder in the direction of $\theta$ and the left-hand side inequality is sharp if and only if $K$ is a double cone in the direction of $\theta$.

\smallskip

\noindent The next proposition is well known. There, 
for
$x,y>0$, 
$ B(x,y) = \int_{0}^{1} \lambda^{x-1} (1-\lambda)^{y-1} d \lambda = \frac{\Gamma(x) \Gamma(y)}{\Gamma(x+y)} $ is the Beta function and $\Gamma(x) = \int_{0}^{\infty} \lambda^{x-1} e^{-\lambda} d\lambda $ is the Gamma function.
\par
\begin{proposition} \label{prop1}
\noindent Let $K$ be a symmetric convex body in $\mathbb R^n$ of volume $1$. Let $ 1\ls p <\infty$ and $\theta \in S^{n-1}$. Then 
$$ 
B(p+1, n)^{\frac{1}{p}} \ls  \frac{h_{Z_{p}(K)} (\theta) }{h_{K}(\theta) } \ls \left( \frac{n}{p+1}\right)^{\frac{1}{p}}  .
$$
\end{proposition}

\noindent {\bf Proof.} As $|K|=1$, 
$$ \frac{2}{n} h_{K}(\theta) f_{K,\theta}(0) \ls 1     \ls  2 h_{K}(\theta) f_{K,\theta}(0) .  $$
\noindent 
Hence, on the one hand, with (\ref{Brunn}), 
\begin{eqnarray*}
h_{Z_{p}(K)}^{p}(\theta) &=& 2  \int_{0}^{h_{K}(\theta) } t^{p}  f_{K,\theta}(t) d t \ls 2  f_{K,\theta}(0) \int_{0}^{h_{K}(\theta) } t^{p} d t \\
& = &
  \frac{2}{p+1}   f_{K,\theta}(0) \  h_{K}^{p+1}(\theta) \ls    \frac{n}{p+1}  h_{K}^{p}(\theta) . 
\end{eqnarray*}
\noindent 
On the other hand, also with with (\ref{Brunn}), 
\begin{eqnarray*}
  h_{Z_{p}(K)}^{p}(\theta) &=&2  \int_{0}^{h_{K}(\theta) } t^{p}  f_{K,\theta}(t) d t \gr 2 f_{K,\theta}(0)\   \int_{0}^{h_{K}(\theta) } t^{p}  \left(1- \frac{t}{h_{K}(\theta)} \right)^{n-1}     d t \\
  & =&
2  f_{K,\theta}(0) h_{K}^{p+1}(\theta) \int_{0}^{1} s^{p} (1-s)^{n-1} d s \gr B(p+1, n) h_{K}^{p}(\theta)  .
\end{eqnarray*}
The proof is complete. $\hfill\Box $ 

\vskip 4mm
As it was mentioned in the introduction, it was proved in  \cite{PW}  that if $K$ is a $C^{2}_{+} $ symmetric convex body of volume $1$, then 
$$ \lim_{p\rightarrow \infty} \frac{p}{\log{p}} \left( |Z_{p}^{\circ}(K) | - |K^{\circ}|\right) = \frac{n(n+1)}{2}|K^{\circ}| . $$
Before we consider the case of  polytopes, we show that for every convex body we have that $ |Z_{p}^{\circ}(K) | - |K^{\circ}| = O (\frac{p}{\log{p}}). $
\noindent 
In particular, the following proposition holds.
\vskip 3mm
\begin{proposition} \label{prop2}
\noindent Let $K$  be a symmetric convex body in $\mathbb R^n$ of volume $1$. Then 
$$ n|K^{\circ}| \ls \lim_{p\rightarrow \infty} \frac{p}{\log{p}} \left( |Z_{p}^{\circ}(K) | - |K^{\circ}|\right)  \ls n^{2}|K^{\circ}| . $$
\end{proposition}
\par 
\noindent {\bf Proof.} We have that 
\begin{eqnarray*}
  |Z_{p}^{\circ}(K) | - |K^{\circ}| &=& \frac{1}{n}  \int_{S^{n-1}} \frac{1}{h_{Z_{p}(K)}^{n} (\theta)}- \frac{1}{h_{K}^{n} (\theta)} d\sigma(\theta) \\
  &= &
 \frac{1}{n}  \int_{S^{n-1}}  \frac{1}{h_{K}^{n} (\theta)} \left( \frac{h_{K}^{n}(\theta)}{h_{Z_{p}(K)}^{n}(\theta)} - 1\right)d\sigma(\theta), 
  \end{eqnarray*}
where $\sigma$ is the usual surface area  measure on $S^{n-1}$.  
By Proposition \ref{prop1}, 
$$   \frac{h_{K}^{n}(\theta)}{h_{Z_{p}(K)}^{n}(\theta)} \gr \left( \frac{n}{p+1}\right)^{-\frac{n}{p}} = 1+ \frac{ n\log{p} }{p} \pm o (\frac{p}{\log{p}})$$
and
$$  \frac{h_{K}^{n}(\theta)}{h_{Z_{p}(K)}^{n}(\theta)}  \ls  B(p+1, n)^{-\frac{n}{p}} = 1+ \frac{ n^{2}\log{p} }{p} \pm o (\frac{p}{\log{p}}). $$
For the last equality see e.g. \cite{PW}, Lemma 4.3 - which is also stated here as Lemma \ref{betafunction}.
Lebesgue's convergence theorem  completes the proof. $\hfill\Box $ 

\smallskip

\section{Polytopes}

\noindent 
Let $K$ be a convex  polytope in $\mathbb R^n$ with vertices $v_1, \dots, v_M$.
For $0 \leq k \leq n-1$, 
let $\mathcal{A}_k= \{F_k: F_k \  \text{is a k-dimensional face of }\  K\}$. 
For $\theta  \in S^{n-1}$ and $ 0 \leq s \leq h_k(\theta)$
let
\begin{equation*}\label {g-func}
g(\theta, s) = \text{card}\left(\{v_i:  v_i \in K \cap \{\langle v_i, \theta \rangle \geq s\}\right).
\end{equation*}
\par
\noindent
Let 
\begin{equation}\label{bad}
{\cal{B}}_K =\{ \theta \in S^{n-1}:  \forall \  s \leq h_K(\theta): g(\theta, s) >1\}
\end{equation}
and
\begin{equation}\label{good}
{\cal{G}}_K =\{ \theta \in S^{n-1}:  \exists \   s <  h_K(\theta): g(\theta, s) =1\}
\end{equation}
Finally, for $\theta \in \mathcal{G}_K$, let
\begin{equation}\label{stheta}
s_\theta = \min \{ s >0:  g(\theta, s) =1\}
\end{equation}

\vskip4mm
\noindent
{\bf Remarks.} Let  $\theta \in  \mathcal {G}_K$.
\vskip 2mm
\noindent
 (i)
Then there is a vertex $v_i$ such that for all $s_\theta \leq s \leq h_K(\theta)$
\begin{equation*}
 \{ x\in K : \langle x, \theta \rangle \geq  s \} = \text{co} \big[ K \cap (\theta^{\perp} + s  \theta ), v_i \big]
 \end{equation*}
\par
\noindent
(ii)  Recall that $f_{K, \theta}(s)= | K \cap (\theta^{\perp} + s  \theta )|$.
We have for all $s_\theta \leq s \leq h_K(\theta)$
\begin{equation}\label{vol-section}
f_{K, \theta}(s) = f_{K, \theta}(s_\theta) \left(\frac{1-\frac{s}{h_K(\theta)}}{1-\frac{s_\theta}{h_K(\theta)}}\right)^{n-1}
\end{equation}

\vskip 4mm
For a convex body $K$,  let $H_K = \max_{\theta \in S^{n-1}} h_K(\theta)$. 
\newline
For $1 \leq k \leq n$,  let $K$ be a $k$-dimensional convex body in a $k$-dimensional affine space of  $\mathbb{R}^n$.
Let 
\begin{equation}\label{inradius}
r(K) = \sup \{   r >0: \exists  \ x \in K \hskip 2mm  \text{such that} \hskip 2mm x+r B^k_2 \subseteq K\}
\end{equation}
be the inradius of $K$.
 Let 
 \begin{equation*}
 r_0 = \min_{1 \leq k \leq n-1} \min_{F_k \in \mathcal{A}_k} r(F_{k}) 
  \end{equation*}
Note that   $r_0 >0$. We also 
put $ h_0 = \max_{u \in \mathcal{B}_K} h_{K}(u) $.
\vskip 2mm
\noindent 
For $ \delta>0$,   we define
\begin{equation}\label{Adelta}
A(\delta)  = \{ \theta \in S^{n-1}:  \exists \  u \in \mathcal{B}_K : \| \theta - u\| < \delta \}  .
\end{equation}
and
\begin{equation}\label{sdelta}
s(\delta) = \sup_{\theta \in S^{n-1} \setminus A(\delta)} \frac{s_\theta}{h_K(\theta)}
\end{equation}

\vskip 2mm
\noindent
{\bf Remark.}
$s(\delta) <1$ and if $\theta \rightarrow \phi$ where $\phi \in \mathcal {B}_K$, then by continuity,   $\frac{s_\theta}{h_K(\theta)} \rightarrow 1$.
Hence we may assume that for  $\delta >0$ small enough, $s(\delta)$ is attained on the ``boundary" of  $S^{n-1} \setminus A(\delta)$.
\vskip 3mm
\noindent
\begin{lemma} \label{lemma1}
\noindent Let $K$ be a $0$-symmetric  polytope in $\mathbb R^n$ of volume $1$.  Then 
for $\delta $ small enough,
$$
s(\delta)=\sup_{\theta \in S^{n-1} \setminus A(\delta)} \frac{s_\theta}{h_K(\theta)}  \leq 1 - \frac{\delta r_0 }{2h_0}
$$
\end{lemma}

\noindent {\bf Proof.} Let $\delta \leq \frac{h_0}{H_K}$. 
By the above Remark, for $\delta >0$ small enough,  there exists  $\phi \in  S^{n-1} \setminus A(\delta)$   such that  $s(\delta)= 
\frac{s_\phi}{h_K(\phi)}$.
\newline
As  $\phi \in  S^{n-1} \setminus A(\delta)$, there exists $u \in \mathcal{B}_K$, such that $\|u - \phi\| = \delta$.
Let 
$v \in \partial K$ be that vertex of $K$ such that $\langle \phi, v \rangle = \max_{x \in K} \langle \phi, x \rangle$.
Let  
$$x_0 = \{ \alpha \phi: \alpha \geq 0\} \cap \partial K, \hskip 5mm z_0 = \{ \alpha u: \alpha \geq 0\} \cap \partial K,
$$ 
and
$$
d_1 = \|x_0 - z_0\|, \hskip 5mmd_2=  \|x_0 - v\|.
$$
$x_0$, $v$ and $z_0$ lie in the $n-1$-dimensional face $F$ orthogonal to $u$.
As $\phi \in \mathcal{G}_K$, we may also assume that $\delta$ is small enough such that  
$ s _\phi = \|x_0\|$, and hence $s(\delta) = \frac{\|x_0\|}{h_K(\phi)}$.
\par
Let $\omega$ be the angle between $ \phi$ and $ u$. Then 
$$ \tan{\omega}= \frac{d_{1}}{h_{K}(u)} \ \ {\rm and} \ \  \sin{\omega}= \frac{h_{K}(\phi) - s_\phi}{d_{2}} . $$
Hence
$$
\frac{h_{K}(\phi) - s_\phi}{d_{2}}  = \frac{d_1 \cos \omega}{h_{K}(u)}
$$
and thus
$$
\frac{s_\phi}{h_{K}(\phi) } = 1- \frac{d_1d_2 \cos\omega}{h_K(u) h_K(\phi)}.
$$
As $d_2 \geq r_0$ and as $\delta \leq \frac{d_1 \cos\omega}{h_K(u)}$,
we get that 
$$
\frac{s_\phi}{h_{K}(\phi) } \leq  1- \frac{\delta \  r_0 }{h_K(\phi)}.
$$
Now observe that 
$$
h_k(\phi) = h_K(\phi -u) + h_K(u) \leq \delta H_K + h_K(u)   \leq 2 h_0.
$$
Therefore,
$$
\frac{s_\phi}{h_{K}(\phi) } \leq  1- \frac{\delta  r_0 }{2h_0}.
$$
$\hfill\Box $ 

\vskip 4mm
Let $f: \mathbb R_{+} \rightarrow \mathbb R_{+} $ be a $C^{2}$ $\log$-concave function with $\int_{\mathbb R_{+}} f(t) dt < \infty$ and let $p\ge 1$. 
Let $g_p(t) = t^p f(t)$ and let $t_p =t_p(f)$ the unique point such that $g^{\prime}(t_p)=0$. We make use of the following Lemma due to B. Klartag \cite{Kl1}  (Lemma 4.3 and Lemma 4.5).

\begin{lemma} \label{klar}
Let $f$ be as above.  For every $\varepsilon \in (0,1)$,
$$ \int_{0}^{\infty} t^p f(t) dt \leq \left( 1+ Ce^{-cp\varepsilon^{2}} \right) \int_{t_p(1-\varepsilon)}^{t_p(1+\varepsilon)} t^p f(t) dt $$
where $C>0$ and  $c>0$ are universal constants.
\end{lemma}

\vskip 4mm
We will use   Lemma \ref{klar} for the function $f_{K, \theta}(s)= | K \cap (\theta^{\perp} + s  \theta )|$ in the proof of the next lemma.
First we oberve
\vskip 3mm
\noindent
\begin{remark}  \label{remark2}
Let $\theta  \in \mathcal {G}_K$. 
As above, let $g_p(t)= t^p f_{K, \theta}(t)$  and let $t_p $ be the unique point such that $g_p^{\prime}(t_p)=0$. Note that, 
since $  t_{p} \rightarrow h_{K}(\theta)$, as $p\rightarrow \infty$ (see e.g. \cite{PW}, Lemma 4.5),    for $p$ large enough - namely $p$ so large that $t_p \geq s_\theta$ -  we can use 
(\ref{vol-section}) and compute $t_p$.
\begin{equation}\label{tp}
t_p = \frac{p}{p+n-1} h_K(\theta)
\end{equation}
\end{remark}

\vskip 3mm
\noindent
We will also use (see e.g.  \cite{PW}, Lemma 4.3).
\vskip 3mm
\noindent
\begin{lemma} \label{betafunction}
\vskip 2mm
\noindent
 Let $p >0$. Then
\begin{eqnarray*}
 \left(B\left(p+1, n\right)\right)^\frac{n}{p} &=&
 1- \frac{n^2}{p} \log p + \frac{n}{p} \log \left(\Gamma(n)\right) +
 \frac{n^4}{2p^2} (\log p )^2 - 
\frac{n^3}{p^2}   \log \left(\Gamma(n)\right) \log p \\
 & \pm& o(p^2).
\end{eqnarray*} 
\end{lemma}

\vskip 3mm

\begin{lemma} \label{lemma2}
\noindent 
Let $K$ be a  $0$-symmetric polytope in $\mathbb R^n$ of volume $1$. For all sufficiently small $\delta$,  for all $\theta \in S^{n-1} \setminus A(\delta)$ and for all  $p \geq \frac{\alpha_n(K)}{\delta}$, we have 
\begin{eqnarray*}
 \left(\frac{ h_{Z_{p}(K)}(\theta) }{h_K(\theta)} \right)^n  \leq 1- n^2\  \frac{ \log p}{p} + (n-1)n\  \frac{ \log\frac{1}{\delta}}{p} + \frac{c_{K,n} }{p}.
\end{eqnarray*}
$\alpha_n(K)=\frac{4(n-1) h_0}{r_0}$ and $c_{K,n}$ are  constants that depend on $K$ and $n$ only.
\end{lemma}

\noindent {\bf Proof.} 
Let $0< \delta \leq \frac{h_0}{H_K}$ be as in Lemma \ref{lemma1}. Let $\theta \in S^{n-1} \setminus  A(\delta)$. Hence, in particular, $\theta \in \mathcal{G}_K$.
By Lemma \ref{klar} we have for all $\varepsilon \in (0,1)$
\begin{eqnarray*}
 h_{Z_{p}(K)}^{p}(\theta) &=& 2  \int_{0}^{h_{K}(\theta) }  t^{p}  f_{K,\theta}(t) d t \\
 &\leq& 2\left( 1+ Ce^{-c p\varepsilon^2}\right) \int_{(1-\varepsilon) t_{p} }^{h_{K}(\theta)} t^{p}  f_{K,\theta}(t) d t
 \end{eqnarray*}
Since $  t_{p} \rightarrow h_{K}(\theta)$, as $p\rightarrow \infty$ (see e.g. \cite{PW}, Lemma 4.5), there exists  $p_{\varepsilon} >0$ (which we will now determine),  
such that  for all $p \geq p_\varepsilon$, 
\begin{equation}\label{pepsilon}
 (1-\varepsilon) t_{p} \geq s_{\theta}.
 \end{equation}  
\par
\noindent
By (\ref{tp}), (\ref{pepsilon}) holds for all $p \geq p_\varepsilon$ with 
$$
p_\varepsilon \geq   \frac{(n-1)\frac{s_\theta}{h_K(\theta} }{1-\varepsilon - \frac{s_\theta}{h_K(\theta} }.
$$  
By Lemma \ref{lemma1},  $\frac{s(\theta)}{h_K(\theta)}  \leq 1 - \frac{\delta r_0 }{2h_0}$
and thus (\ref{pepsilon}) holds for all $p \geq p_\varepsilon$ with 
$$
p_\varepsilon \geq  \frac{n-1}{\delta}\  \frac{2h_0 - \delta r_0 }{r_0  - 2 h_0 \varepsilon/\delta}.
$$
We choose $\varepsilon = \frac{r_0 \delta }{4h_0}$. Then for 
$$
p_\varepsilon \geq  \frac{n-1}{\delta}\  \frac{4h_0}{r_0}
$$
the estimate (\ref{pepsilon}) holds for all $p \geq p_\varepsilon$ uniformly for all $\theta \in S^{n-1} \setminus A(\delta)$.
Thus, using also (\ref{vol-section}), 
\begin{eqnarray}\label{above1}
 h_{Z_{p}(K)}^{p}(\theta) &\leq&    2\left( 1+ Ce^{-c p\varepsilon^2}\right) \int_{(1-\varepsilon) t_{p} }^{h_{K}(\theta)} t^{p}  f_{K,\theta}(t) d t \nonumber \\
& \leq &
  2 \left( 1+ Ce^{-cp\varepsilon^2}\right) \int_{s_\theta}^{h_{K}(\theta)} t^{p}  f_{K,\theta}(t) d t \nonumber \\
&=&2 \left( 1+ Ce^{-cp \varepsilon^2}\right) \  \frac{ h_{K}^{p+1}(\theta)  f_{K,\theta} (s_{\theta}) }{\left(1-\frac{s_\theta}{h_K(\theta)}\right)^{n-1}} \  \int_{\frac{s_{\theta}}{h_K(\theta)}}^{1}  u^{p} \left( 1-u \right)^{n-1}d u \nonumber  \\
&\leq&2 \left( 1+ Ce^{-cp \varepsilon^2}\right) \  \frac{ h_{K}^{p+1}(\theta)  f_{K,\theta} (0) }{\left(1-\frac{s_\theta}{h_K(\theta)}\right)^{n-1}} \  \int_{\frac{s_{\theta}}{h_K(\theta)}}^{1}  u^{p} \left( 1-u \right)^{n-1}d u \nonumber \\
& \leq &
n \  \left( 1+ Ce^{-c p\varepsilon^2 }\right) \ B(p+1, n) \ h_{K}^{p}(\theta) \left(\frac{2h_0}{\delta r_
0  }\right)^{n-1}. 
\end{eqnarray}
In the last inequality we have used that $1-\frac{s_\theta}{h_K(\theta)} \geq \frac{\delta r_0 }{2h_0}$ and that $\frac{2}{n} h_{K}(\theta)  f_{K,\theta} (0) \leq |K| = 1$. Equivalently, (\ref{above1}) becomes
\begin{eqnarray*}\label{above2}
\left(\frac{ h_{Z_{p}(K)}(\theta) }{h_K(\theta)} \right)^n \leq  n^\frac{n}{p}  \left( 1+ Ce^{-c p\varepsilon^2 }\right) ^\frac{n}{p} \left( \frac{2h_0}{\delta r_0  }\right)^\frac{(n-1)n}{p} B(p+1, n)^\frac{n}{p}.
\end{eqnarray*}
With  Lemma \ref{betafunction},  we then get
\begin{eqnarray*}
\left(\frac{ h_{Z_{p}(K)}(\theta) }{h_K(\theta)} \right)^n \leq 1- n^2 \frac{ \log p}{p} + (n-1)n \frac{  \log \frac{1}{\delta}}{p} + \frac{c_{K,n} }{p}. 
\end{eqnarray*}
$\hfill\Box $

\vskip 5mm
Let $ \delta \in [0, 1)$ and $ \theta \in S^{n-1}$. We define the cap $C(\theta, \delta)$ of the sphere $S^{n-1}$ around $\theta$  by
$$ C(\theta, \delta):= \{ \phi\in S^{n-1}: \|\phi -\theta \|_{2} \ls \delta \} . $$
We will  estimate the surface area of a cap, and to do so 
we will make use of the following fact which follows immediately from e.g. Lemma 1.3 in \cite{SW1}.
\vskip 4mm
\begin{lemma} \label{area-cap}
Let $\theta \in S^{n-1}$ and $ \delta <1$.
Then
\begin{eqnarray*}
&&\hskip -10mm \text{vol}_{n-1}\big(B^{n-1}_2\big)  \left(1- \frac{\delta^{2}}{4}\right )^{\frac{n-1}{2}} \  \delta^{n-1} \ls \\
&& \hskip 36mm \sigma(C(\theta, \delta)) \leq \\
&&\hskip 40mm \text{vol}_{n-1}\big(B^{n-1}_2\big)  \left(1- \frac{\delta^{2}}{4}\right )^{\frac{n-1}{2}}  \frac{\left(1+ \frac{\delta^{4}}{4}\right )^{\frac{1}{2}} }{\left(1- \frac{\delta^{2}}{2}\right)} \  \delta^{n-1}.
\end{eqnarray*}
\end{lemma}

\noindent 
\vskip 5mm

\noindent
 {\bf Proof of Theorem \ref{theorem1}.}
 \par
 \noindent
For $p$ given, let $\delta= \frac{1}{\log p}$. Let $A(\delta)$ as defined in (2.10). Let $p_{0}$ be  such that $p_0$ and $\delta= \frac{1}{\log p}$ 
satisfy the assumptions of Lemma \ref{lemma2}, 
i.e. $\frac{p_0}{\log p_0} \geq \frac{4(n-1)h_0}{r_0}$. 
By Lemma \ref{lemma2}, we have for all $p \geq p_0$, 
\begin{eqnarray*}
&& |Z_{p}^{\circ}(K) | - |K^{\circ}| 
 \geq \frac{1}{n} \int_{S^{n-1}\setminus A(\delta)} \frac{1}{h_{Z_{p}(K)}^{n} (\theta)} \left( 1- \frac{h_{Z_{p}(K)}^{n}(\theta)}{h_{K}^{n}(\theta)}\right) d\sigma(\theta) \\
& &\geq   \frac{1}{n} \int_{S^{n-1}\setminus A(\delta)}\frac{1}{h_{Z_{p}(K)}^{n} (\theta)} \bigg(
 \frac{n^{2} \log{p}}{p} - (n-1)n \frac{\log \log p}{p} + \frac{c_{K,n} }{p}
\bigg) d\sigma(\theta) \\
  &&=  \frac{1}{n} \int_{S^{n-1}}\frac{1}{h_{Z_{p}(K)}^{n} (\theta)} \bigg(
 \frac{n^{2} \log{p}}{p} - (n-1)n \frac{\log \log p}{p} + \frac{c_{K,n} }{p} 
\bigg)d\sigma(\theta)\\
&&-  \frac{1}{n} \int_{A(\delta)}\frac{1}{h_{Z_{p}(K)}^{n} (\theta)} \bigg(
 \frac{n^{2} \log{p}}{p} - (n-1)n \frac{\log \log p}{p} + \frac{c_{K,n} }{p} 
\bigg) d\sigma(\theta).
\end{eqnarray*}
Hence,
\begin{eqnarray*}
&& \frac{p}{\log p} \left(|Z_{p}^{\circ}(K) | - |K^{\circ}|\right) 
\geq \\
&& \frac{1}{n} \int_{S^{n-1}}\frac{1}{h_{Z_{p}(K)}^{n} (\theta)} \bigg(
 n^{2}- \frac{(n-1)n \log \log p}{\log p} + \frac{c_{K,n} }{\log p} 
\bigg) d\sigma(\theta)\\
&&-  \frac{1}{n} \int_{A(\delta)}\frac{1}{h_{Z_{p}(K)}^{n} (\theta)} \bigg(
 n^{2}- \frac{(n-1)n \log \log p}{\log p} + \frac{c_{K,n} }{\log p}  
\bigg)d\sigma(\theta).
\end{eqnarray*}
\par
\noindent
Note that, since $K$ is centrally symmetric, $r(K)= \inf_{\theta \in S^{n-1}} h_K(\theta)$. Also, since $Z_p(K)$ converges to $K$,
for $p$ sufficiently large, $h_{Z_{p}(K)}^{n} (\theta) \geq \left( \frac{r(K)}{2}\right)^n$ for every $\theta \in S^{n-1}$.
Together with  Lemma \ref{area-cap} we thus get
\begin{eqnarray*}
&&\frac{1}{n} \int_{A(\delta)}\frac{1}{h_{Z_{p}(K)}^{n} (\theta)} d\sigma(\theta) \leq \\
&&\frac{2^{n+1}}{n\  r(K)^n} \text{card} \left(\mathcal{B}_K\right) \text{vol}_{n-1} \left(B^{n-1}_2\right) \delta^{n-1}
\left(1- \frac{\delta^{2}}{4}\right )^{\frac{n-1}{2}}  \frac{\left(1+ \frac{\delta^{4}}{4}\right )^{\frac{1}{2}} }{\left(1- \frac{\delta^{2}}{2}\right)} \\
&& \leq \frac{2^{n+1} \text{card} \left(\mathcal{B}_K\right) }{n\  r(K)^n}
\frac{\text{vol}_{n-1} \left(B^{n-1}_2\right) }{(\log p)^{n-1}}.
\end{eqnarray*}
 
By Proposition \ref{prop2} and  Lebesgue's convergence theorem we can interchange integration and limit
and get
\begin{eqnarray*}
&&\lim_{p\rightarrow \infty} \frac{p}{\log{p}} \left( |Z_{p}^{\circ}(K) | - |K^{\circ}|\right)   \geq \\
&&  \frac{1}{n} \int_{S^{n-1}} \lim_{p\rightarrow \infty} \frac{1}{h_{Z_{p}(K)}^{n} (\theta)} \bigg(
 n^{2}- \frac{(n-1)n \log \log p}{\log p} + \frac{c_{K,n}}{\log p}  
\bigg) d\sigma(\theta) \\
&&-  \frac{2^{n+1}  \text{card} \left(\mathcal{B}_K\right) \text{vol}_{n-1} \left(B^{n-1}_2\right)}{n\  r(K)^n} 
  \lim_{p\rightarrow \infty}  
\bigg( \frac{n^{2}}{(\log p)^{n-1}}- \frac{(n-1)n \log \log p}{(\log p)^{n}} + \frac{c_{K,n}}{(\log p)^{n}} 
\bigg)\\
&&= n^{2}|K^{\circ}| . 
\end{eqnarray*}
Here, we have also used that $\lim _{ p \rightarrow \infty} h_{Z_{p}(K)} (\theta) = h_K(\theta)$.
\vskip 3mm
\noindent The inequality from above follows by Proposition \ref{prop2}.
 $\hfill\Box $ 

\section{Approximation with uniformly convex bodies}

\noindent Let $K$ be a symmetric convex body in $\mathbb R^n$ and $2\ls p < \infty$. We say that $ K$ is $p$-uniformly convex (with constant $C_{p}$) (see e.g. \cite{Garl, LindTz}), if for every $x, y \in \partial K$, 
$$ \| \frac{ x+y}{2}\|_{K} \ls 1 - C _{p}\|x-y\|_{K}^{p} . $$
\noindent We will need the following Proposition. The proof is based on Clarkson inequalities and can be found in e.g. (\cite{Garl}, pp. 148).

\begin{proposition}\label{uniformconvex}
\noindent Let $K$ be a compact set in $\mathbb R^n$ of volume $1$. Then for $p\gr 2$, $Z_{p}^{\circ}(K)$ is $p$-uniformly convex with constant $C_{p}= \frac{1}{p2^{p}}$. 
\end{proposition}

\smallskip

\noindent The symmetric difference metric between two convex bodies $K$ and $C$ is 
$$ d_{s}(C,K) = | ( C\setminus K) \cup (K\setminus C)| . $$

\smallskip

\noindent 
{\bf Proof of Theorem \ref{theorem2}.} 
\par
\noindent
Let $P_{1} = \frac{P^{\circ}} {|P^{\circ}|^{\frac{1}{n}}} $. Then $P_{1}^{\circ} = |P^{\circ}|^{\frac{1}{n}} P$ and $ |P_{1}^{\circ}|= |P| |P^{\circ}| $. Let $K_{p} = |P^{\circ}|^{-\frac{1}{n}} Z_{p}^{\circ}(P_{1})$. Then by Proposition \ref{uniformconvex} we have that $K_{p}$ is uniformly convex. Note that $P \subseteq K_{p}$. By Theorem 1.1 we have that 
$$ \lim_{p\rightarrow \infty} \frac{ p}{\log{p}} \left( |Z_{p}^{\circ} (P_{1})| - | P_{1}^{\circ}| \right) = n^{2} | P_{1}^{\circ}| . $$ 
So, for every $\varepsilon>0$, there exists $p_{0}( \varepsilon, P)$ such that 
$$  d_{s}(P, K_{p}) = |K_{p}| - |P| = \frac{ 1}{| P^{\circ}| } \left( |Z_{p}^{\circ}(P_{1})| - |P_{1}^{\circ}| \right) \ls  $$
$$ (1+ \varepsilon) n^{2} \frac{|P_{1}^{\circ}|}{ |P^{\circ}|} \frac{\log{p}}{p}=  (1+ \varepsilon) n^{2} |P|  \frac{\log{p}}{p} .$$
We choose $\varepsilon =1 $ and the proof is complete. $\hfill\Box $

\vskip 2mm 
\noindent 
Grigoris Paouris\\
{\small Department of Mathematics}\\
{\small Texas A \& M University}\\
{\small College Station, TX, , U. S. A.}\\
{\small \tt }\\ \\
\noindent
\and 
Elisabeth Werner\\
{\small Department of Mathematics \ \ \ \ \ \ \ \ \ \ \ \ \ \ \ \ \ \ \ Universit\'{e} de Lille 1}\\
{\small Case Western Reserve University \ \ \ \ \ \ \ \ \ \ \ \ \ UFR de Math\'{e}matique }\\
{\small Cleveland, Ohio 44106, U. S. A. \ \ \ \ \ \ \ \ \ \ \ \ \ \ \ 59655 Villeneuve d'Ascq, France}\\
{\small \tt elisabeth.werner@case.edu}\\ \\


\begin{thebibliography}{00}

\bibitem{Gar} {\sc R.J.~Gardner}, {\em Geometric tomography}, Encyclopedia of
Mathematics and its Applications, {\bf 58}. Cambridge University
Press, Cambridge, 1995.



\bibitem{GaZ}
{\sc R. J. Gardner and G. Zhang},
{\em Affine inequalities and radial mean bodies.}
 Amer. J. Math. {\bf 120} no.3 (1998), 505-528.

\bibitem{Garl} {\rm D.\ J.\ H.\ Garling}, {\sl Inequalities, A journey into Linear Analysis}, Cambridge University Press  (2007). 

\bibitem{GrZh}
{\sc
E. Grinberg and G. Zhang,} {\em Convolutions, transforms, and convex bodies}, Proc. 
London Math. Soc. {\bf 78} no.3 (1999), 77-115.


\bibitem{Hab}
{\sc C. Haberl}, {\em Blaschke valuations}, Amer. J. of Math., in press



\bibitem{HabSch}
{\sc C. Haberl and F.  Schuster,} {\em General Lp affine isoperimetric inequalities}.
J. Differential Geometry {\bf 83} (2009), 1-26.

\bibitem{HLYZ}
{\sc C. Haberl, E. Lutwak, D. Yang and G. Zhang,} {\em The even Orlicz Minkowski problem}, 
Adv. Math. {\bf 224} (2010), 2485-2510




\bibitem{Klain1}
{\sc D. Klain},
{\em  Star valuations and dual mixed volumes}, Adv. Math. {\bf 121} (1996), 80-101. 

\bibitem{Klain2}
{\sc D. Klain},
{\em Invariant valuations on star-shaped sets},   Adv. Math.  {\bf 125} (1997), 95-113. 

\bibitem{Kl1} {\rm B. Klartag}, {\sl A central limit theorem for convex
sets,} Invent. Math. {\bf 168} (2007) 91-131.

\bibitem{LindTz}
{\sc J. Lindenstrauss and L. Tzafriri}, {\em  Classical Banach spaces I, II},  reprint, Springer, 1996.



\bibitem{Lud}
{\sc M. Ludwig}, {\em  Minkowski valuations},  
Trans. Amer. Math. Soc.  {\bf 357}  (2005), 4191-4213.



\bibitem{Lud3}
{\sc M. Ludwig}, {\em
Minkowski areas and valuations}, 
J. Differential Geometry, {\bf 86} (2010), 133-162.

\bibitem{LR1}
{\sc M. Ludwig and M. Reitzner,}  {\em A Characterization of Affine Surface Area}, Adv.  Math. {\bf 147} (1999), 138-172.

\bibitem{LR2}
{\sc M. Ludwig and M. Reitzner,} {\em A classification of $SL(n)$
invariant valuations.}  Ann. of Math. {\bf 172 } (2010), 1223-1271. 


\bibitem{Lu2}{\sc E. Lutwak}, {\em The Brunn-Minkowski-Firey theory II : Affine and
geominimal surface areas}, Adv. Math. {\bf 118}  (1996),   244-294.

\bibitem{LZ} {\rm E. Lutwak and G. Zhang}, {\sl Blaschke-Santal\'{o}
inequalities}, J. Differential Geom. {\bf 47} (1997), 1--16.


\bibitem{MW} {\rm M. Meyer and E. Werner}, {\sl The Santalo-regions of a convex body},
Transactions of the AMS {\bf 350}   (1998), 4569--4591.
 


\bibitem{PW} {\rm G. Paouris and E. Werner}, {\sl Relative entropy of cone measures and $L_p$-centroid bodies},
Proc. London Math. Soc., in press (available at arXiv:0909.436).

\bibitem{Sch} {\rm R. Schneider}, {\sl Convex Bodies: The
Brunn-Minkowski Theory}, Encyclopedia of Mathematics and its
Applications {\bf 44}, Cambridge University Press, Cambridge (1993).

\bibitem{Schm} {\rm M. Schmuckenschl\"ager}, {\sl The distribution function of the convolution
square of a convex symmetric body in $\mathbb R^n$}, Israel Journal of Math. {\bf 78}
(1992), 309--334.

\bibitem{SW2} {\rm C. Sch\"utt and E. Werner}, {\sl The convex floating body},   
{Math. Scand.} {\bf 66}, (1990), 275--290.



\bibitem{SW1} {\rm C. Sch\"utt and E. Werner}, {\sl Polytopes with Vertices Chosen Randomly
from the Boundary of a Convex Body},   Geom. Aspects of Funct. Analysis,
Lecture Notes in Math. {\bf 1807} (2003), 241--422.


\bibitem{SW2004}
{\sc C. Sch{\"u}tt and E. Werner}, {\em Surface bodies and
p-affine surface area.} Adv. Math. {\bf 187}  (2004), 98-145.

\bibitem{SA1}
{\sc A. Stancu,} {\em The Discrete Planar $L_0$-Minkowski
Problem.} Adv. Math. {\bf 167} (2002),  160-174.

\bibitem{SA2}
{\sc A. Stancu}, {\em On the number of solutions to the
discrete two-dimensional $L_0$-Minkowski problem.} Adv. Math. {\bf
180} (2003), 290-323.


 
 \bibitem{W1} {\rm E. Werner}, {\sl  Illumination bodies and affine surface area}, Studia
Math. {\bf 110}  (1994), 257--269 .
 

\bibitem{WY} {\sc E. Werner and D. Ye}, {\em New $L_{p}$ affine isoperimetric inequalities},
Adv. Math. {\bf 218} no.3 (2008), 762-780.


\bibitem{WY1} {\sc E. Werner and  D. Ye}, {\em Inequalities for mixed $p$-affine surface area},
{Math. Ann.} {\bf  347}  (2010), 703-737

\bibitem{Z3}{\sc G. Zhang}, {\em New Affine  Isoperimetric Inequalities}, ICCM 2007, Vol. II, 239-267.

\medskip
\end{thebibliography}
\end{document}